\documentclass[12pt]{article}

\makeatletter
\newfont{\footsc}{cmcsc10 at 8truept}
\newfont{\footbf}{cmbx10 at 8truept}
\newfont{\footrm}{cmr10 at 10truept}
\renewcommand{\ps@plain}{%
\renewcommand{\@oddfoot}{\footsc the electronic journal of combinatorics {\footbf 11} (2004), \#R00\hfil\footrm\thepage}}
\makeatother
\usepackage{amssymb}
\usepackage{graphicx}
\usepackage{amsmath}
\usepackage{amsthm}
\usepackage{array}
\usepackage{comment}

\renewcommand{\P}{\mathbb{P}}
\newcommand{\Q}{\mathbb{Q}}

\newtheorem{theorem}{Theorem}[section]

\newtheorem{lemma}[theorem]{Lemma}
\newtheorem{corollary}[theorem]{Corollary}
\newtheorem{conjecture}[theorem]{Conjecture}

\theoremstyle{definition}
 
\newtheorem{example}{Example}[section]

\theoremstyle{remark}

\numberwithin{equation}{section}

\renewcommand{\ge}{\geqslant}
\renewcommand{\geq}{\geqslant}
\renewcommand{\le}{\leqslant}
\renewcommand{\leq}{\leqslant}

\begin{document}
\renewcommand{\theequation}{\arabic{section}.\arabic{equation}}

\thispagestyle{empty}

\vskip 20pt
\begin{center}
{\bf BOTTOM SCHUR FUNCTIONS}
\vskip 15pt
{\bf Peter Clifford
} \\
{\it CNRI,}
{\it Dublin Institute of Technology,
Ireland}\\
{\texttt{peterc@alum.mit.edu}}
\\

\hskip 0.5cm

{\bf Richard P. Stanley \footnote{Partially supported by NSF grant
\#DMS-9988459.}}\\
{\it Department of Mathematics, Massachusetts Institute of
Technology}\\
{\it Cambridge, MA 02139, USA}\\
{\texttt{rstan@math.mit.edu}}\\[.2in]
{\small Submitted: Nov 19, 2003; Accepted 27 Aug, 2004\\MR Subject Classifications: 05E05, 05E10}\\
\end{center}

\begin{abstract}
We give a basis for the space spanned by the sum $\hat{s}_\lambda$ of
the lowest degree terms in the expansion of the Schur symmetric
functions $s_\lambda$ in terms of the power sum symmetric functions
$p_\mu$, where deg$(p_i)=1$. These lowest degree terms 
correspond to minimal border strip tableaux of $\lambda$.  The
dimension of the space spanned by $\hat{s}_\lambda$, where $\lambda$
is a partition of $n$, is equal to the number of partitions of $n$
into parts differing by at least 2. Applying    
the Rogers-Ramanujan identity, the generating function also counts the
number of partitions of $n$ into parts $5k+1$ and $5k-1$.

We also show that a symmetric function closely related to
$\hat{s}_\lambda$ has the same coefficients when expanded in terms of
power sums or augmented monomial symmetric functions. 
\end{abstract}

\section{Introduction}

Let $\lambda = (\lambda_1, \lambda_2, \ldots)$ be a partition of the
integer $n$, i.e., $ \lambda_1 \ge \lambda_2 \ge \cdots \ge 0$ and
$\sum \lambda_i =n$.   The
\emph{length} $\ell(\lambda)$ of a partition $\lambda$ is the 
number of
nonzero parts of $\lambda$. The (Durfee or Frobenius) \emph{rank} of
$\lambda$, denoted rank($\lambda$), is the length of the main diagonal
of the diagram of $\lambda$, or equivalently, the largest integer $i$
for which $\lambda_i \ge i$.  The rank of $\lambda$ is the least
integer  $r$ such that
$\lambda$ is a disjoint union of $r$ border strips (defined below).

Nazarov and Tarasov \cite[Sect. 1]{Na98}, in connection with tensor
products of Yangian modules, defined a generalization of rank to skew
partitions (or skew diagrams) $\lambda / \mu$.  The paper
\cite[Proposition 2.2]{St02} gives several simple equivalent
definitions of rank($\lambda / \mu$).  One of the definitions is that
rank($\lambda / \mu$) is the least integer $r$ such that $\lambda /
\mu$ is a disjoint union of $r$ border strips. It develops a general
theory of minimal border strip tableaux of skew shapes, introducing
the concepts of the snake sequence and the interval set of a skew
shape $\lambda / \mu$.  These tools are used to count the number of
minimal border strip decompositions and minimal border strip tableaux
of $\lambda / \mu$.  In particular, the paper \cite{St02} gives an
explicit combinatorial formula for the coefficients of the $p_{\nu}$,
where $\ell(\nu) = \mbox{rank} (\lambda / \mu)$, which appear in the
expansion of $s_{\lambda / \mu}$.

The paper \cite{St02} considered a degree operator
deg$(p_{\nu})=\ell(\nu)$ and defined the \emph{bottom Schur functions}
to be the sum of the terms of lowest degree which appear in the
expansion of $s_{\lambda / \mu}$ as a linear combination of the
$p_{\nu}$.  We study the bottom Schur functions in detail when $\mu =
\emptyset$.  In particular, in Section \ref{sec:spacebotschur} we give
a basis for the vector space they span.

In Section \ref{sec:identity} we show that when we substitute $ip_i$
for $p_i$ in the expansion of a bottom Schur function in terms of
power sums, then the resulting symmetric function has the same
coefficients when expanded in terms of power sums or augmented
monomial symmetric functions. 
\section{Definitions}\label{sec:definitions}

\begin{figure}[h]
\centering
\includegraphics[width=2.0in]{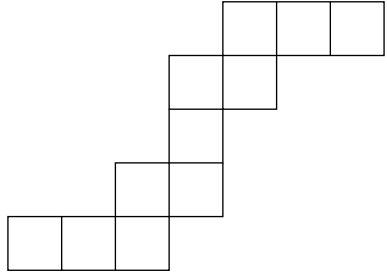}
\caption{The border strip $75443/4332$}
\label{fig:borderstrip}
\end{figure}

In general we follow \cite[Ch.~7]{St99} for notation and terminology
involving symmetric functions. Let
$\lambda$ be a partition of $n$ with Frobenius rank $k$.
Recall that $k$ is the length of the main diagonal
of the diagram of $\lambda$, or equivalently, the largest integer $i$
for which $\lambda_i \ge i$.  Let
$m_i(\lambda)  = \# \{j: \lambda_j = i \}$, the
number of parts of $\lambda$ equal to $i$.  Define $z_\lambda =
1^{m_1(\lambda)} m_1(\lambda)! 2^{m_2(\lambda)} 
  m_2(\lambda)!\cdots$.    A \emph{border strip} (or \emph{rim hook}
  or \emph{ribbon}) is a connected skew shape with no $2 \times 2$
  square.   An example is $75443/4332$ whose diagram is illustrated in
  Figure \ref{fig:borderstrip}.  Define the \emph{height} ht($B$) of a
  border strip $B$ to be one less than its number of rows.

\begin{figure}[h]
\centering
\includegraphics[width=1.5in]{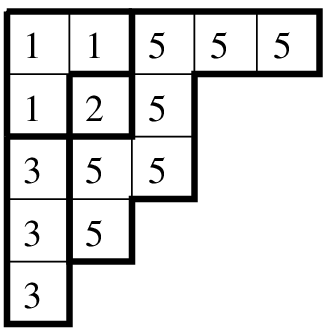}
\caption{A border strip tableau of $53321$ of type $(3,1,3,0,7)$}
\label{fig:borderstriptableau}
\end{figure}

Let $\alpha=(\alpha_1, \alpha_2, \ldots)$ be a weak composition of
$n$, i.e., $\alpha_i\geq 0$ and $\sum \alpha_i=n$.
Define a \emph{border strip tableau} of shape $\lambda$ and type
$\alpha$ to be an assignment of positive integers to the squares of
$\lambda$ such that:
\begin{enumerate}
\item every row and column is weakly increasing,
\item the integer $i$ appears $\alpha_i$ times, and
\item the set of squares occupied by $i$ forms a border strip.
\end{enumerate}
Equivalently, one may think of a border-strip tableau as a sequence
$\emptyset=\lambda^0 \subseteq \lambda^1 \subseteq \cdots \subseteq
\lambda^r \subseteq \lambda$ of partitions such that each skew shape
$\lambda^i / \lambda^{i+1}$ is a border-strip of size $\alpha_i$.
For instance,  Figure \ref{fig:borderstriptableau} shows a border
strip tableau of $53321$ of type $(3,1,3,0,7)$.
 It is easy to
see (in this nonskew case) that the smallest number of strips in a
border-strip tableau is rank($\lambda$).  Define the \emph{height}
ht($T$) of a border-strip tableau $T$ to be  
$$\mbox{ht}(T) = \mbox{ht}(B_1)+\mbox{ht}(B_2)+ \cdots +
\mbox{ht}(B_k)$$
where $B_1, \ldots, B_k$ are the (nonempty) border strips appearing in
$T$.  In the example we have ht$(T)= 1+0+2+3 = 6$.  Now we can define
$$\chi^{\lambda}(\nu) = \sum_T (-1)^{\mathrm{ht}(T)},$$
summed over all border-strip tableaux of shape $\lambda$ and type
$\nu$.  Since there are at least rank$(\lambda)$ strips in every
tableau, we have that $\chi^{\lambda}(\nu) = 0$ if $\ell(\nu) <
$ rank$(\lambda)$.  The numbers $\chi^{\lambda} (\nu)$ for $\lambda,
\nu \vdash n$ are the values of the irreducible characters
$\chi^\lambda$ of the symmetric group $S_n$.

Finally we can express 
the Schur function $s_\lambda$ in terms of power sums $p_{\nu}$,
viz., 
  \begin{equation} s_\lambda = \sum_{\nu} \chi^{\lambda}
     (\nu) \frac{p_{\nu}}{z_{\nu}}. \label{eq:char} 
  \end{equation}

Define deg($p_i$) $=1$, so deg($p_{\nu}) = \ell(\nu)$.  The
\emph{bottom Schur function} $\hat{s}_\lambda$ is defined to be the
lowest degree part of $s_\lambda$, so  
 $$ \hat{s}_\lambda = \sum_{\nu: \ell(\nu) = \mathrm{rank}(\lambda)}
  \chi^{\lambda} (\nu) \frac{p_{\nu}}{ z_{\nu}}. $$
Also write $\tilde{p}_i = \frac{p_i}{i}$.  For instance,
$$s_{321} = \frac{1}{45}p_1^6 -\frac{1}{9}p_3p_1^3+
\frac{1}{5}p_1p_5- \frac{1}{9}p_3^2.
$$
Hence
\begin{eqnarray*}
\hat{s}_{321} & = & \frac{1}{5}p_1p_5- \frac{1}{9}p_3^2 \\
            & = & \tilde{p}_1 \tilde{p}_5- \tilde{p}_3^2.
\end{eqnarray*}

\begin{figure}[h]
\centering
\includegraphics[width=2.0in]{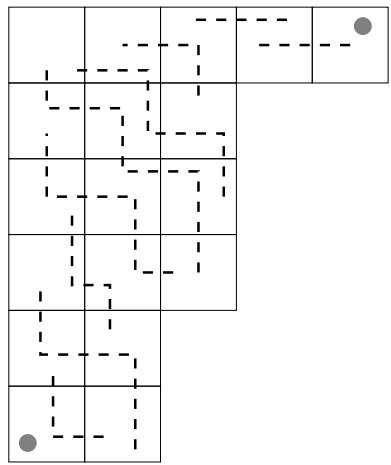}
\caption{Snakes for the shape 533322}
\label{fig:snakes}
\end{figure}

We identify a partition $\lambda$ with its \emph{diagram}
  $$ \lambda=\{(i,j)\,:\,1\leq j\leq \lambda_i\}. $$
Let $e$ be an edge of the lower envelope of $\lambda$, i.e., no square
of $\lambda$ has $e$ as its upper or left-hand edge.  We will define a
certain subset $S_e$ of squares of $\lambda$, called a \emph{snake}.
If $e$ is horizontal and $(i,j)$ is the square of $\lambda$ having $e$
as its lower edge, define
\begin{eqnarray}
 S_e & = & (\lambda) \cap \{ (i,j), (i-1,j), (i-1, j-1), \nonumber \\
     &   & (i-2, j-1),(i-2,j-2),\ldots \label{eq:left} \}.
\end{eqnarray}
If $e$ is vertical and $(i,j)$ is the square of $\lambda$ having $e$
as its right-hand edge, define
\begin{eqnarray}
 S_e & = & (\lambda) \cap \{ (i,j), (i,j-1), (i-1, j-1), \nonumber \\
     &   & (i-1, j-2),(i-2,j-2),\ldots \label{eq:right} \}.
\end{eqnarray}
In Figure \ref{fig:snakes} the nonempty snakes of the shape $533322$
are shown with dashed paths through their squares, with a single
bullet in the two snakes with just one square.  The \emph{length}
$\ell(S)$ of a snake $S$ is one fewer than its number of squares; a
snake of length $i-1$ (so with $i$ squares) is call an \emph{i-snake}.
Call a snake of \textbf{even} length a \emph{left snake} if $e$ is horizontal
and a \emph{right snake} if $e$ is vertical.
It is clear that the snakes are linearly ordered from lower left to
upper right.  In this linear ordering, replace a left snake with the
symbol $L$, a right snake with $R$, and a snake of odd length with
$O$.  The resulting sequence (which does not determine $\lambda$)
is called the \emph{snake sequence} of $\lambda$, denoted
SS$(\lambda)$.  
For instance, from Figure \ref{fig:snakes} we see that
$$ \mbox{SS}(533322) = LLOOLORROOR.$$

\begin{lemma} \label{lemma2.1}
The $L$'s in the snake sequence correspond exactly to horizontal edges
of the lower envelope of $\lambda$ which are below the line $x+y=0$.  The
$R$'s correspond 
exactly to vertical edges of the lower envelope of $\lambda$ which are
above the line $x+y=0$.  All other edges of the lower envelope of
$\lambda$ are labelled by $O$'s.
\end{lemma}
Clearly we could have defined the snake sequence this way;
however, the definitions above also hold for skew shapes.
Lemma~\ref{lemma2.1} only holds when $\lambda$ is a straight (i.e.,
nonskew) shape.
\begin{proof}
Let $e$ be an edge of the lower envelope of $\lambda$ below the line
$x+y =0$.  Let $(i,j)$ be the square of $\lambda$ having $e$ as its
lower edge.  The last square in the snake is some square in the first
column of $\lambda$.  So if $e$ is horizontal then the last square is
$(i-j+1,1)$, the snake has an odd number of squares and so has even
length, and is labelled by $L$.  If $e$ is vertical then the last
square is $(i-j,1)$, the snake has an even number of squares, so has
odd length, and is labelled by $R$.  The case when $e$ is above
$x+y=0$ is proved similarly.
\end{proof}

\begin{corollary}\label{cor:lsfirst}
In the snake sequence of $\lambda$, the $L$'s occur strictly to the
left of the $R$'s.  
\end{corollary}

The number of horizontal edges of the lower envelope of $\lambda$
which are below the line $x + y=0$ equals the length of the main
diagonal of the diagram of $\lambda$, which is the rank of $\lambda$.
Similarly the number of vertical edges of the lower envelope of
$\lambda$ which are above the line $x+y=0$ also equals the rank of
$\lambda$.  Henceforth we fix $k=\mathrm{rank}(\lambda)$.

\begin{figure}[h]
\centering
\includegraphics[width=3.0in]{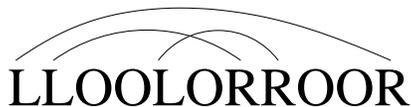}
\caption{An interval set of the shape 533322}
\label{fig:iset}
\end{figure}

Let SS$(\lambda) = q_1 q_2 \cdots q_m$, and define an \emph{interval
  set} of $\lambda$ to be a collection ${\cal I}$ of $k$ ordered
pairs,
$$ {\cal I} = \{(u_1,v_1),\ldots,  (u_k,v_k) \} ,$$
satisfying the following conditions:
\begin{enumerate}
\item the $u_i$'s and $v_i$'s are all distinct integers,
\item $1 \le u_i < v_i \le m$,
\item $q_{u_i} = L$ and $q_{v_i} = R$.
\end{enumerate}
Figure \ref{fig:iset} illustrates the interval set
$\{(1,11),(2,7),(5,8) \}$ of the 
shape 533322.

Given an interval set
${\cal I} = \{ (u_1,v_1) , \ldots ,  (u_k,v_k) \} $,
define the \emph{crossing number} $c({\cal I})$
to be the number 
of \emph{crossings} of ${\cal I}$, i.e. the number of pairs $(i,j)$
for which $u_i<u_j<v_i<v_j$.
 
Let {\boldmath $T$}  be a border strip tableau of
shape $\lambda$.  Recall that
$$ \mbox{ht} ( \mbox{\boldmath $T$}  ) = \sum_{B} \mbox{ht}(B),$$
where $B$ ranges over all border strips in {\boldmath $T$} and
$\mbox{ht}(B)$ is one less than the number of rows of $B$.  Define
$z(\lambda)$ to be the height $\mbox{ht}  ( \mbox{\boldmath $T$}  ) $
of a ``greedy border strip tableau'' {\boldmath $T$} of shape
$\lambda$ obtained by starting with $\lambda$ and successively
removing the largest possible border strip.  (Although  {\boldmath
  $T$} may not be unique, the set of border strips appearing in
{\boldmath $T$} is unique, so $\mbox{ht}  ( \mbox{\boldmath $T$}  ) $
is well-defined.)

The connection between bottom Schur functions and interval sets was
given by Stanley \cite[Theorem 5.2]{St02}: 
$$ \hat{s}_{\nu} = (-1)^{z(\nu)} \sum_{ {\cal I} = \{(u_1,v_1),\ldots,
  (u_k,v_k) \} }  (-1)^{c( {\cal I} )} \prod_{i=1}^{k}
\tilde{p}_{v_i-u_i}, $$
where ${\cal I}$ ranges over all interval sets of $\nu$.
 
For example the shape $321$ has snake sequence $LOLROR$.  There are
two interval sets, $\{(1,4),(3,6)\}$ with crossing number $1$, and
$\{(1,6),(3,4)\}$ with crossing number $0$.   So  as we saw before
$$
\hat{s}_{321}  =  \tilde{p}_1 \tilde{p}_5- \tilde{p}_3^2.
$$

\section{Bottom Schur Functions of straight shapes}\label{sec:botschur}

\begin{lemma}\label{l:shapeord} The lexicographic order on shapes $\nu$
  whose length $\ell(\nu)$ equals their rank $k$ is
  equal to the reverse lexicographical order (with respect to the
  ordering L$<$R$<$O) on their snake sequences.
\end{lemma}
\begin{proof}
Since $\ell(\nu)=k$, the snake sequence begins with $k$ L's.  If 
 the length of the $i$th row
of $\nu$ is $k+j$, then there are $j$ O's to the left of the $(k-i+1)$st R.
\end{proof}
Denote the \emph{complete homogeneous symmetric functions} by 
$h_\lambda$.
Recall that the Jacobi-Trudi  identity expresses the $s_\lambda$'s in terms of the
$h_\mu$'s: 
 $$ s_\lambda= \mbox{det}(h_{\lambda_i - i +j})_{i,j=1}^n,$$
where we define $h_i =0$ for $i < 0$. 
For example 
$$s_{554421} =  \mbox{det} \left[
\begin{array}{llllll}
h_5 & h_6 & h_7 & h_8 & h_9 & h_{10} \\
h_4 & h_5 & h_6 & h_7 & h_8 & h_9 \\
h_2 & h_3 & h_4 & h_5 & h_6 & h_7 \\
h_1 & h_2 & h_3 & h_4 & h_5 & h_6 \\
0   &  0  & 1   & h_1 & h_2 & h_3 \\
0   &  0  & 0   &  0  &  1 & h_1
\end{array}
\right].
$$

Since $h_n = \sum_{\lambda \vdash n} \frac{p_\lambda}{z_\lambda}$,
the term of lowest degree (in $p$) in the expansion of a given $h_n$
in terms of the $p_j$ is just $\frac{p_n}{n}=\tilde{p}_n$.
For a product $h_{n_1} h_{n_2} \cdots h_{n_j}$ the term of lowest
degree in the expansion in terms of the $p_j$ is just $\tilde{p}_{n_1}
\tilde{p}_{n_2} 
\cdots \tilde{p}_{n_j}$.  So we have that 
$\hat{s}_\lambda=$ terms of lowest order in $\mbox{det}
  (\tilde{p}_{\lambda_i   - i
  +j})_{i ,j=1}^n$ (since
the $p_\lambda$ are algebraically
independent, and since $\mbox{det}(h_{\lambda_i - i +j}) = s_\lambda \ne
0$, this determinant will not vanish).
For example
$$\hat{s}_{554421} =  \mbox{ terms of lowest order in det} \left[
\begin{array}{llllll}
\tilde{p}_5 & \tilde{p}_6 & \tilde{p}_7 & \tilde{p}_8 & \tilde{p}_9 &
\tilde{p}_{10} \\ 
\tilde{p}_4 & \tilde{p}_5 & \tilde{p}_6 & \tilde{p}_7 & \tilde{p}_8 &
\tilde{p}_9 \\ 
\tilde{p}_2 & \tilde{p}_3 & \tilde{p}_4 & \tilde{p}_5 & \tilde{p}_6 &
\tilde{p}_7 \\ 
\tilde{p}_1 & \tilde{p}_2 & \tilde{p}_3 & \tilde{p}_4 & \tilde{p}_5 &
\tilde{p}_6 \\ 
0   &  0  & 1   & \tilde{p}_1 & \tilde{p}_2 & \tilde{p}_3 \\
0   &  0  & 0   &  0  &  1 & \tilde{p}_1
\end{array}
\right].
$$
Since $p_0 = 1$, the terms of lowest order are those which contain the
most  number of 1's.

Row $i$ of the matrix will have a $1$ in position $(i,j)$ if
$\lambda_i-i+j=0$, i.e. if $\lambda_i<i$ (this shows that the number
of rows of $JT_\lambda$ which do not contain a $1$ is another
definition of $\mbox{rank}(\lambda)$ \cite[Prop.~2.2]{St02}).

Let $JT^*_p$ be the
matrix obtained from the original Jacobi-Trudi matrix by removing
every row and column which contains a 1 and replacing the $h_i$ with
$\tilde{p}_i$.  We show below that this matrix is not singular and so
we have 
$$\hat{s}_\lambda= \mbox{det } JT^*_p.$$
For example
$$\hat{s}_{554421} =  \mbox{det} \left[
\begin{array}{llll}
\tilde{p}_5 & \tilde{p}_6 & \tilde{p}_8 & \tilde{p}_{10} \\ 
\tilde{p}_4 & \tilde{p}_5 & \tilde{p}_7 & \tilde{p}_9 \\ 
\tilde{p}_2 & \tilde{p}_3 & \tilde{p}_5 & \tilde{p}_7 \\ 
\tilde{p}_1 & \tilde{p}_2 & \tilde{p}_4 & \tilde{p}_6 \\ 
\end{array}
\right].
$$

Any minor of the Jacobi-Trudi matrix for a
shape  $\lambda$ is the
Jacobi-Trudi matrix for some skew shape $\mu / \sigma$.  For let $JT^*$ be
some minor of size $m$ of some Jacobi-Trudi matrix $JT$.  If the
entry in position $(i,j)$ is $h_{x}$ put $jt^*_{i,j} = x$.
Now we can set 
$$\sigma_{i} = jt^*_{1,m}-jt^*_{1,i}-m+i,$$
and
$$\mu_i = jt^*_{i,i} + \sigma_i.$$
Again note that since the $p_\lambda$ are algebraically independent
and $\mbox{det } JT^* = s_{\mu / \sigma} \ne 0$, we have $ \mbox{det }
JT^*_p \ne 0$.  

In our running example, we have $\sigma_1 = 10-5-4+1 = 2, \sigma_2 =
10 -6 - 4 +2 = 2, \sigma_3 = 10 -8 -4 + 3 = 1 $ and $ \sigma_4 = 10 -10
-4+4 = 0$.  Hence $\sigma = (2210)$.  Also $\mu_1 = 5+2, \mu_2 = 5+2, \mu_3
= 5+1 $ and $ \mu_4= 6+0$.  Thus $\mu=(7766)$.   Therefore we have that
$\hat{s}_{554421}$ equals the determinant of the Jacobi-Trudi matrix
of $7766 / 2210$ with the $h$'s replaced by $\tilde{p}$'s.  

\begin{lemma}\label{l:squarek}
If the skew shape $\mu / \sigma$ has the Jacobi-Trudi matrix $JT^*$
obtained by removing all rows and columns with a 1 from a
Jacobi-Trudi matrix $JT$ of a shape $\lambda$ with rank $k$, then $\mu /
\sigma$ contains a square of size $k$.
\end{lemma}
\newpage 
The rank of $554421$ is $4$, and the diagram of $7766 / 2210$ 
does indeed contain a square of size 4:
\begin{figure}[h]
\centering
\includegraphics[width=1.0in]{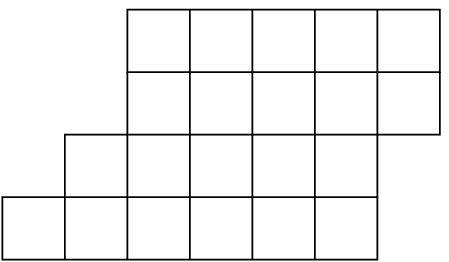}
\end{figure}
\begin{proof}
We give a proof due to Christine Bessenrodt, 
greatly improving our  
original proof.  Define 
$\mu'_i=\ell(\lambda)-k+\lambda_i \ \ (i=1,\ldots,k)$ 
and $ \sigma'_i=\#\{s|\lambda_s \leq k-i\} \ \ (i=1,\ldots,k)$.  
We give a diagrammatic definition of $\mu'$ and $\sigma'$ 
which also illustrates that the 
skew diagram $\mu'/\sigma'$ contains a square of size $k$. 
Consider $\lambda$ as a $k \times k$ square with 
two partitions $\alpha$ and $\beta$ glued to it, i.e. $\lambda=(k+\beta_1,
\ldots, k+\beta_k, \alpha_1, \ldots, \alpha_{\ell(\lambda)-k})$.  Flip 
$\alpha$ over its anti-diagonal and then 
glue the bottom right corner of the result to the bottom left corner 
of the square.  
The final diagram is the skew diagram of $\mu'/\sigma'$.
We show that 
$\mu = \mu'$ and $\sigma = \sigma'$.

The $k$ rows of $JT^*$ are contained in the first $k$ rows of $JT$, so
$\mu_i = \lambda_i +c$ for some constant 
$c$.  The last column of $JT$ does not have
a 1 in it, so it will not be removed,
and its first $k$ entries will be the last column of $JT^*$.  Hence
$jt^*_{1,k} = jt_{1,\ell(\lambda)} = \lambda_1 + \ell(\lambda)-1$.  Since
$jt^*_{1,k} = \mu_1-\sigma_k + k -1$, we have 
$\mu_i = \ell(\lambda)-k+\lambda_i=\mu'_i$.

The first $k$ entries of the last column of $JT$ are retained.  
Then we remove the next
$\#\{s| \lambda_s = 1\}$
columns to its left, do not remove the next column, remove the next
$\#\{s| \lambda_s = 2\}$
columns to the left,
and so on.  
Formally we have 
$jt^*_{1,k-j} = jt^*_{1,k-j+1}-1-\#\{s|\lambda_s=j\}$. 
Combining this with $\sigma_{i} = jt^*_{1,k}-jt^*_{1,i}-k+i$
gives us 
$\sigma_i = \#\{s|\lambda_s \leq k-i\} = \sigma'_i$.

\end{proof}

\section{The space spanned by the bottom Schur functions}
\label{sec:spacebotschur}
Before we use the above results to give a basis for the
space spanned by the bottom Schur functions, we must first recall some
classical tableaux theory.

If $\lambda/\mu$ is a skew shape, then a \emph{standard Young tableau}
(SYT) of shape $\lambda/\mu$ is a labelling of the squares of
$\lambda/\mu$ with the numbers $1,2, \ldots, n$, each number appearing
once, so that every row and column is increasing.  A
\emph{semistandard Young tableau} (SSYT) of shape $\lambda/\mu$ is a
labelling of the squares of $\lambda/\mu$ with positive integers that
is weakly increasing in every row and strictly increasing in every
column.  We say that $T$ has \emph{type} $\alpha = (\alpha_1,
\alpha_2, \ldots)$ if $T$ has $\alpha_i$ parts equal to $i$.
\begin{figure}[h]
\centering
\includegraphics[width=2.0in]{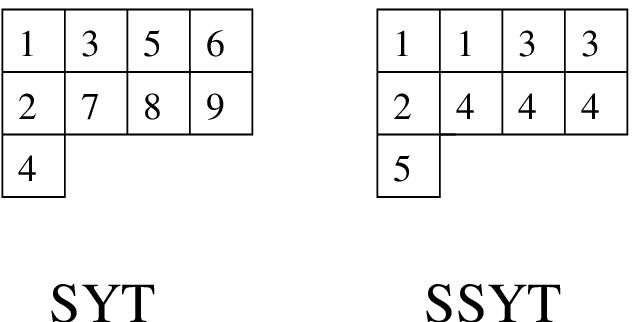}
\end{figure}

\begin{figure}[h]
\centering
\includegraphics[width=2.25in]{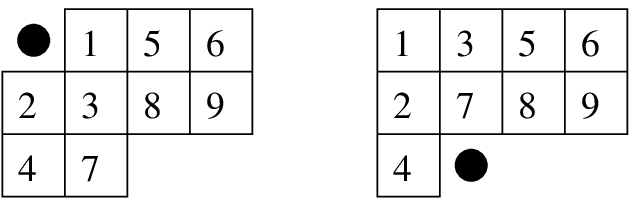}
\caption{Jeu de taquin slides}
\label{fig:jeudetaquin}
\end{figure}

Now we define an operation (of Sch\"utzenberger) on standard Young
tab\-leaux called a \emph{jeu de taquin slide}.  Given a skew shape
$\lambda / \mu$, consider the squares $b_0$ that can be added to
$\lambda / \mu$, so that $b_0$ shares at least one edge with $\lambda
/ \mu$, and $\{b_0 \} \cup \lambda / \mu$ is a valid skew shape.
Suppose that $b_0$ shares a lower or right edge with $\lambda / \mu$
(the other situation is completely analogous).  There is at least one
square $b_1$ in $\lambda / \mu$ that is adjacent to $b_0$; if there
are two such squares, then let $b_1$ be the one with a smaller entry.
Move the entry occupying $b_1$ into $b_0$.  Then repeat this
procedure, starting at $b_1$.  The resulting tableau will be a
standard Young tableau.  Analogously if $b_0$ shares an upper or left
edge, the operation is the same except we let $b_1$ be the square with
the bigger entry from two possibilities.  For example we illustrate
both situations in Figure \ref{fig:jeudetaquin}; the tableau on the
right results from playing jeu de taquin beginning at the square
marked by a bullet on the tableau on the left (and vice versa).

Two tableaux $T$ and $T'$ are called \emph{jeu de taquin equivalent} if
one can be obtained from another by a sequence of jeu de taquin
slides. Given an SYT $T$ of shape $\lambda/\mu$, there is exactly one
SYT $P$ of straight shape, denoted jdt$(T)$, that is jeu de taquin
equivalent to $T$ \cite[Thm.~A1.2.4]{St99}.

The \emph{reading word} of a
(semi)standard Young tableau is the sequence of entries of $T$
obtained by concatenating the rows of $T$ bottom to top.  For example,
the tableau on the left in Figure \ref{fig:jeudetaquin}
has the reading word $472389156$.
 The \emph{reverse reading word} of a
tableau is simply the reading 
word read backwards.

A \emph{lattice permutation} is a sequence $a_1 a_2 \cdots a_n$ such
that in any initial factor $a_1 a_2 \cdots a_j$, the number of $i$'s
is at least as great as the number of $i+1$'s (for all $i$).  For
example $123112213$ is a lattice permutation.

The \emph{Littlewood-Richardson coefficients} $c^{\lambda}_{\mu \nu}$
are the coefficients in the expansion of a skew Schur function in the
basis of Schur functions:
$$ s_{\lambda / \mu} = \sum_{\nu} c^{\lambda}_{\mu \nu} s_{\nu} .$$
The Littlewood-Richardson rule is a combinatorial description of the
coefficients  $c^{\lambda}_{\mu \nu}$.  We will use two different
versions of the rule.

\begin{theorem}[Sch\"utzenberger, Thomas]
 Fix an SYT $P$ of shape $\nu$.  The
  Littlewood-Richardson coefficient $c^{\lambda}_{\mu \nu}$ is equal
  to the number of SYT of shape $\lambda / \mu$ that are jeu de taquin
  equivalent to $P$.
\end{theorem} 
\begin{figure}[h]
\centering
\includegraphics[width=1.25in]{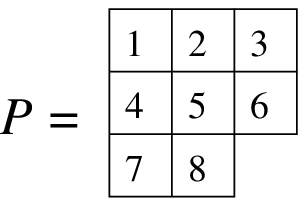}
\end{figure}
For example, let $\lambda=(5,3,3,1), \mu = (3,1), $ and $ \nu =
(3,3,2)$.  Consider the tableau $P$ of shape $\nu$ shown above.  There
are exactly two SYTs $T$ of shape $\lambda / 
\mu$ such that $\mathrm{jdt}(T)=P$, namely,
\begin{figure}[h]
\centering
\includegraphics[width=3.5in]{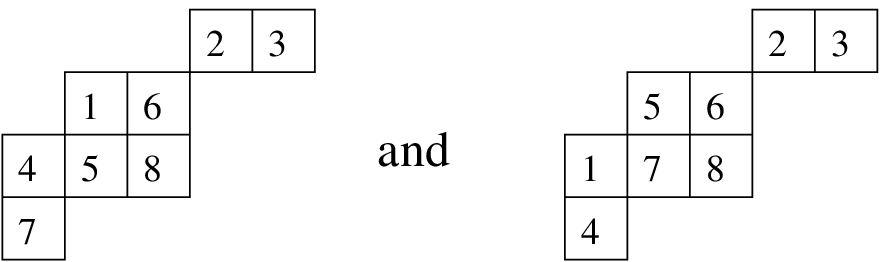}
\end{figure}

\begin{theorem}[Sch\"utzenberger, Thomas]
  The
  Littlewood-Richardson coefficient $c^{\lambda}_{\mu \nu}$ is equal
  to the number of semistandard Young tableaux of shape $\lambda /
  \mu$ and type $\nu$ whose reverse reading word is a lattice
  permutation.  
\end{theorem} 
For example, with $\lambda=(5,3,3,1), \mu = (3,1), $ and $ \nu =
(3,3,2)$ as above, there are exactly two SSYTs $T$ of shape $\lambda /
\mu$ and type $\nu$ whose reverse reading word is a lattice
  permutation:
\begin{figure}[h]
\centering
\includegraphics[width=3.5in]{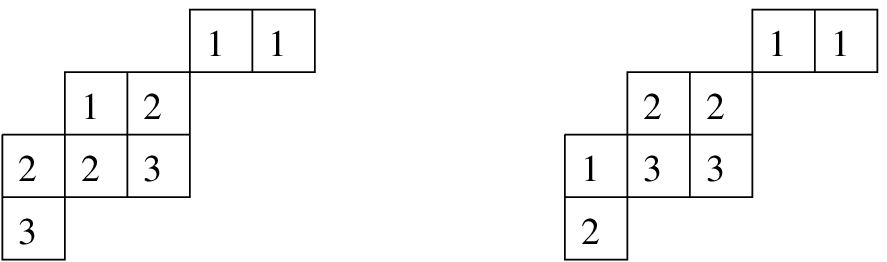}
\end{figure}

Now we have enough machinery to state and prove this section's main theorem.

\begin{theorem}
Fix $n$ and $k$.   The set $ \{ \hat{s}_{\nu} :
\nu \vdash n, \ \mathrm{rank}(\nu) =k $ and $ \ell(\nu) = k \} $ is  
a basis for
the space $\mathrm{span}_{\Q}  \{ \hat{s}_\lambda : \lambda \vdash n \ and  \
\mathrm{rank}(\lambda)  = k \}$.
\end{theorem}
For example if $n=12$ and $k=3$, we have that  $ \{ \hat{s}_{633},
\hat{s}_{543},\hat{s}_{444}\}$ is a basis for $\mathrm{span}_{\Q}  \{
\hat{s}_{633}, \hat{s}_{543}, \hat{s}_{5331}, \hat{s}_{444},
\hat{s}_{4431}, \hat{s}_{4332}, \hat{s}_{43311},
\hat{s}_{3333}, \hat{s}_{33321} , \hat{s}_{333111} \}$.
\begin{proof}
First we prove that the $\hat{s}_{\nu} $ are linearly
independent.  We show that given any such $\nu$, there is
some term in the expansion of $\hat{s}_{\nu} $ which does not
occur in the expansion of any $\hat{s}_{\nu'} $ for $\nu'$
lexicographically less than $\nu$.

From \cite[Theorem 5.2]{St02} we have that 
$$ \hat{s}_{\nu} = (-1)^{z(\nu)} \sum_{ {\cal I} = \{(u_1,v_1),\ldots,
  (u_k,v_k) \} }  (-1)^{c( {\cal I} )} \prod_{i=1}^{k}
\tilde{p}_{v_i-u_i}, $$
where ${\cal I}$ ranges over all interval sets of $\nu$.  Let $t =
\pm p_{j_1 \ge \cdots \ge j_k}$ 
be the term corresponding to the noncrossing interval set ${\cal I}$
  of the 
snake sequence of $\nu$.  We claim that $t$ does not occur in the
expansion of any $\hat{s}_{\nu'} $ for $\nu'$
lexicographically less than $\nu$.  Assume by way of contradiction
that it does occur for some such $\nu'$ with corresponding
interval set ${\cal I}'$.  

Assume inductively that
the first $i-1$ L's are matched with the last $i-1$ R's without
crossings in ${\cal I}'$.  
Let $r_j$ (and $r'_j$ respectively)  be the
position of the $j$th R in the snake sequence of $\nu$ ($\nu'$
respectively).    
By Lemma 
\ref{l:shapeord} $r_j \ge r'_j$.  But the length of the interval
matching the $i$th R from the right in ${\cal I}$ is $r_{k-i+1}-i$.
So for there to be an interval of this length in ${\cal I}'$ we must
match the $i$th R from the right with the $i$th L; this
interval has no crossing.   Proceeding by induction we see that ${\cal
  I}'$ is also noncrossing, and so must equal ${\cal I}$.
This shows that the snake sequences corresponding to $\nu$ and
$\nu'$ are equal.  Identical snake sequences and equal lengths guarantee 
that
$\nu=\nu'$, a contradiction.

Now we prove that the $\hat{s}_{\nu} $ span the space of all
$\hat{s}_\lambda $.   
We have shown
that $\hat{s}_\lambda= \hat{s}_{\mu / \sigma}$.  Expand
$s_{\mu / \sigma}
$ in terms of (straight) Schur functions using the Littlewood Richardson
rule
$$s_{\mu / \sigma} = \sum_{\nu} c_{\sigma \nu}^{\mu}s_{\nu} \ .$$
We need to show that $c_{\sigma \nu}^{\mu}=0$ unless
$\nu$ is of rank $k$ and length $k$.

Fix an SYT P of shape $\nu$.  The Littlewood-Richardson coefficient
$c_{\sigma  \nu}^{\mu}$ is equal to the number of SYT of shape $\mu /
\sigma$ that are jeu de taquin equivalent to P.  Playing jeu de taquin on
a straight-shape tableau of shape $\nu$ can only increase the
length of the shape.  Hence if $c_{\sigma \nu}^{\mu} \ne 0, \ell(\nu) \le
\ell(\mu)=k$. 

The Littlewood-Richardson coefficient
$c_{\sigma  \nu}^{\mu}$ is also equal to the number of semistandard
Young tableaux of shape $\mu / \sigma$ and type $\nu$ whose reverse
reading word is a lattice permutation.  But we know that $\mu / \sigma$
contains a square of size $k$ (by Lemma \ref{l:squarek}).  Therefore the
bottom $k$ boxes of this square must have labels at least $k$.  Since $\ell(\nu)
\le k$, the labels are exactly $k$.  So
$\nu_k \ge k$, i.e. rank($\nu) \ge k$.  Since $\ell(\nu) \le
k$, we must have rank($\nu) = \ell(\nu) = k$.

Taking terms of lowest degree on both sides of
$$s_{\mu / \sigma} = \sum_{\nu} c_{\sigma \nu}^{\mu}s_{\nu} \ ,$$
we have that
$$\hat{s}_\lambda = \hat{s}_{\mu / \sigma} = \sum_{\nu} c_{\sigma
  \nu}^{\mu}\hat{s}_{\nu}  \ ,$$
where the sum is over $\nu$ of length $k$ and rank $k$ as required.
\end{proof}

\section{Dimension of the space spanned by the bottom Schur functions}
\label{sec:dimbotschur}

Let $p_{\le k}(n)$ be the number of partitions of $n$ with length at
most $k$, and define $p_{\le k}(0)=1$.  A partition $\nu \vdash n$ of length 
$k$ and
rank $k$  decomposes into a $k \times k$ square of boxes and a 
partition of $n-k^2$ of length
at most $k$.

\begin{corollary} \label{cor:botdim}
The dimension of the space of bottom Schur functions \\
$\mathrm{span}_{\Q} \{\hat{s}_\lambda:\lambda \vdash n \}$ is 
$$\sum_{k=1}^{\lfloor \sqrt{n} \rfloor} p_{\le k}(n-k^2).$$
\end{corollary}
For example, the first 27 terms in this sequence (beginning with
$n=1$) are
 $$ 1,1,1,2, 2,3,3,4,5,6,7,9, 10,12,14,17,19,23, 26,31,35,41,46,54,61,
    70,79. $$
There is a nice bijection between the above partitions and the set of
partitions $\{ \lambda \vdash n: \lambda_i - \lambda_{i+1} \ge 2 \}$.
For, given a $k$ and a partition $\lambda^* \vdash n-k^2$ with fewer
than $k$ parts, we can set $\lambda_i = \lambda^*_{i} + 2k -2i +1$.
This gives a partition of $n$ with $k$ rows with $\lambda_i - \lambda_{i+1}
\ge 2$ as required.  This is clearly a bijection.  

This classical sequence also gives the number of partitions
of $n$ into parts congruent to $1$ or $4 \mbox{ mod } 5$; equivalently these numbers are the
coefficients 
in the expansion of
       the Rogers-Ramanujan identity
$$
1 + \sum_{n \ge 1}
              \frac{t^{n^2}}{(1-t)(1-t^2) \cdots (1-t^n)} = \prod_{n \ge 1}
              \frac{1}{(1-t^{5n-1})(1-t^{5n-4})}
$$

\section{2-bottom Schur functions}\label{sec:2botschur}
We have shown that a basis for the space spanned by the bottom Schur
functions consists of the $\hat{s}_\lambda$ where $\ell(\lambda) \le
\mbox{rank}(\lambda)$. It is natural to define for fixed $j\geq 1$ the
$j$-\emph{bottom Schur function} $\hat{s}^j_\lambda$ to be the sum of
those terms of degree at most rank$(\lambda)+j-1$ in the expansion
(\ref{eq:char}) (with $\deg p_i=1$ as usual). When $j=2$ we have
verified (using Stembridge's SF package for Maple \cite{St01}) that
for $n\leq 14$ the dimension of the space spanned by
$\{\hat{s}^2_\lambda\,:\,\lambda\vdash n\}$ equals the number of
$\lambda\vdash n$ satisfying
$\ell(\lambda)\leq\mathrm{rank}(\lambda)+1$. This suggests the
following conjecture.

\begin{conjecture}   
A basis for the space spanned by the 2-bottom Schur functions
consists of all 2-bottom Schur functions $\hat{s}^2_\lambda$, where
$\lambda$ is a partition of n satisfying $\ell(\lambda) \leq
\mbox{rank}(\lambda)+1$.
\end{conjecture}

However in the $j=3$ case, the dimensions of the spaces spanned by the
3-bottom Schur functions are $1,2,3,4,6,9,11,15,19,24,30, \ldots$.
We have computed that the numbers of $\lambda\vdash n$ satisfying $\ell(\lambda)
\leq \mbox{rank}(\lambda)+2$ are given by $1, 2, 3, 4,5, 8,
10, 14, 17, 22, 27,
\ldots$. Unfortunately these sequences do not agree.

\section{A condition satisfied by bottom Schur
   functions}\label{sec:identity} 
%
We prove a surprising identity satisfied by a variant of the bottom
Schur functions related to their expansion in terms of power sum and
monomial symmetric functions.

\begin{figure}[h]
\centering
\includegraphics[width=3.0in]{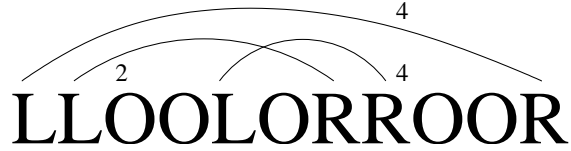}
\caption{A labelled interval set of the shape 533322}
\label{fig:intervalset}
\end{figure}

Fix a shape $\lambda$ of rank $k$.  Given an interval set ${\cal I} =
\{ (u_1,v_1), \ldots, (u_k,v_k) \},$ of $\lambda$
and a labelling of the intervals $(\alpha _1, \ldots, \alpha _k)$ such
that $
\alpha_i \in \P$, 
define 
$$x^{ {\cal I}} = \prod_{i=1}^{k} x_{\alpha_i}^{v_i-u_i}.$$
Recall that $c ( {\cal I} )$ is the number of crossings of the
interval set ${\cal I}$.  Figure
\ref{fig:intervalset} shows a labelled interval set of the shape
533322 with the snake sequence $LLOOLORROOR$.  For this interval set
$c ( {\cal I} ) = 1$ and for this labelling $x^{ {\cal I}} = x_4^{10}
x_2^5 x_4^3 = x_2^5 x_4^{13}$.

\begin{figure}[h]
\centering
\includegraphics[width=5.0in]{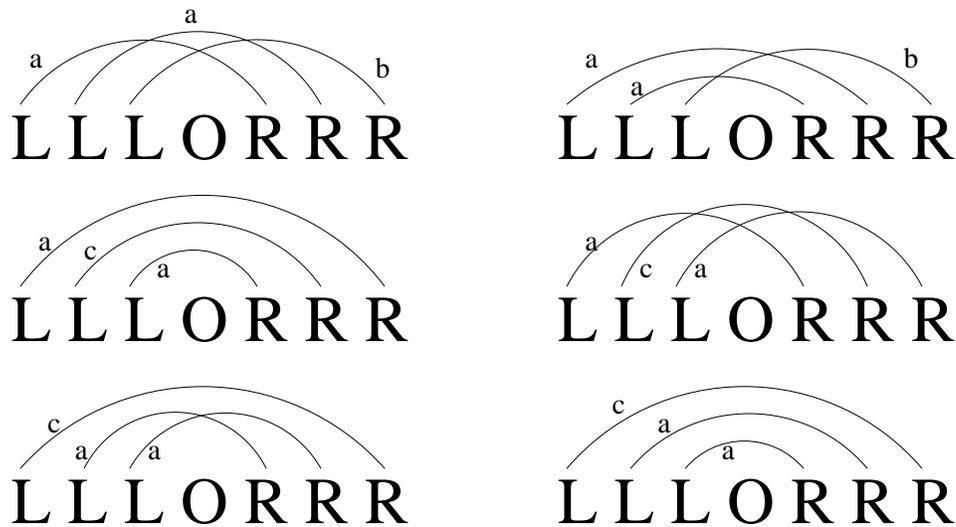}
\caption{Some labelled interval sets of the shape 444}
\label{fig:involution}
\end{figure}

\begin{example}
For the shape $\lambda = (4,4,4)$ with snake sequence $LLLORRR$,
Figure \ref{fig:involution} depicts some of the labelled interval sets.  In the
top left we have $(-1)^{c( {\cal I})}x^{ {\cal I}} = (-1)^3 x_a^4 x_a^4 x_b^4 =
-x_a^8 x_b^4 $.
In the top right  we have $(-1)^{c( {\cal I})}x^{ {\cal I}} = (-1)^2 x_a^5
x_a^3 x_b^4 = x_a^8 x_b^4$. In fact in 
every row the
term  $(-1)^{c( {\cal I})}x^{ {\cal I}}$ in the left column is exactly
the negative    
of the corresponding term $(-1)^{c( {\cal I})}x^{ {\cal I}}$ in the
right column.   
\end{example}

\begin{lemma}\label{l:indistinct}
Fix a shape $\lambda$.  Then
$\sum (-1)^{c( {\cal I})} x^{ {\cal I}} = 0$, where the sum is over
all labelled interval sets of $\lambda$ with at least one label repeated. 
\end{lemma}

\begin{proof}
We give a sign reversing involution on these labelled interval sets.
Examine a specific labelled interval set ${\cal I}$.  
Since we are dealing
with straight (non-skew) shapes,  we know by Corollary \ref{cor:lsfirst}
that the snake 
sequence has all the L's 
before any of the R's, or $u_k < v_1$.  So any two intervals $i$ and
$j$ ($>i$ say) either intersect ($u_i < u_j < v_i < v_j $) or are
nested ($u_i<u_j < v_j < v_i$).  

Let $a$ be the smallest label which is repeated.  The intervals in ${\cal I}$
are ordered by where they start, so identifying the first two
intervals $i$ and $j$ ($>i$) labelled by $a$ is well-defined.

Simply change the interval $(u_i,v_i)$ to $(u_i,v_j)$ and the interval
$(u_j,v_j)$ to $(u_j, v_i)$, while preserving the label $a$ on both.  Where
the intervals start remains unchanged, so these intervals remain the
first two intervals labelled by $a$.  Hence this operation is an
involution.  Note that if the two intervals initially nested, they now
intersect, and if they initially intersected, they now nest.  The parity of 
the number of crossings of these two with any other interval is preserved
under this operation.  So the parity of the total number of crossings has
changed, and this
involution is sign reversing.  The rows of Figure \ref{fig:involution}
are some examples of this involution (if $a \ne c$).

Thus given any labelled interval set with repeated labels, there is a
unique labelled interval set with one more (or fewer) crossings, and
so the sum of all such terms $(-1)^{c( {\cal I})} x^{ {\cal I}}$ is zero.
\end{proof}

Fix a shape $\lambda$.  Recall that 
$$ \hat{s}_\lambda = (-1)^{z(\lambda)} \sum_{ {\cal I} = \{(u_1,v_1),\ldots,
  (u_k,v_k) \} }  (-1)^{c( {\cal I} )} \prod_{i=1}^{k}
\tilde{p}_{v_i-u_i}, $$
where ${\cal I}$ ranges over all interval sets of $\lambda$.  

For
another shape $\mu$, define $c_\mu$ to be the coefficient of
$\tilde{p}_\mu$ in the above sum, i.e.
$$ c_\mu = (-1)^{z(\lambda)} \sum_{ {\cal I}}  (-1)^{c( {\cal I} )},
$$
where the sum is over all interval sets (of $\lambda$) of type $\mu$.  

\begin{lemma}\label{l:powersums}
$c_\mu p_\mu  = (-1)^{z(\lambda)} \sum (-1)^{c( {\cal I})} x^{
  {\cal I}}$, where the sum is over all labelled interval sets of type
$\mu$.  
\end{lemma}

\begin{figure}[h]
\centering
\includegraphics[width=2in]{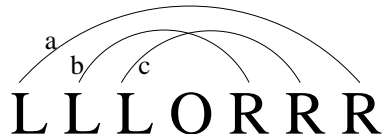}
\caption{A labelled interval set of the shape 444}
\label{fig:labiset}
\end{figure}

\begin{example}\label{e:powersums}  Consider as before the shape $\lambda =
  (4,4,4)$ with snake 
  sequence $LLLORRR$.  In particular, consider the interval
  set $\{ (1,7),(2,5),(3,6) \}$ of type $(6,3,3)$ labelled by
  $(a,b,c)$.  This interval set is illustrated in Figure
  \ref{fig:labiset}.  If $a=b=c$ then $x^{{\cal I} } = x_a^{12}$ and
  the sum over all such labellings will give $x_1^{12} + x_2^{12} +
  \cdots = m_{(12)}$.  If $a=b \ne c$ the sum over all such labellings
  will give $ x_1^6 x_1^3 x_2^3 + x_1^6 x_1^3 x_3^3 + x_2^6 x_2^3
  x_1^3 + \cdots = x_1^9 x_2^3 +x_1^9 x_3^3 + x_2^9 x_1^3 + \cdots =
  m_{93}$.  Similarly if $a=c \ne b$ we will get 
  $m_{93}$, and if $b=c \ne a$ we will get $x_1^6 x_2^3 x_2^3 +
  x_2^6 x_1^3 x_1^3 + \cdots = 2x_1^6 x_2^6 + \cdots =  2m_{66}$.
  Finally if the three labels are all different, the sum will give $x_1^6
  x_2^3 x_3^3 + x_1^6 x_3^3 x_2^3 + \cdots = 2m_{633}$.  So the
  sum over all such labellings is
  $m_{(12)}+2m_{66}+2m_{93}+2m_{633} = p_{633}$.  
\end{example}

\begin{proof}
We need to show that for every interval set ${\cal I}$ of type $\mu$,
$\sum x^{{\cal I} } = p_\mu$,  where the sum is over all labellings of
${\cal I}$.  First note that the intervals can be ordered largest
first and left to right among intervals of the same length.  So the
$i$th interval is well defined, and has length $\mu_i$.

By definition $p_\mu = p_{\mu_1} p_{\mu_2} \cdots
p_{\mu_{\ell(\mu)}} = (x_1^{\mu_1}+x_2^{\mu_1}+
\cdots)(x_1^{\mu_2}+x_2^{\mu_2}+ \cdots) \cdots
(x_1^{\mu_{\ell(\mu)}}+ x_2^{\mu_{\ell(\mu)}}+ \cdots)$.  But if we
expand this product into monomials $x_{i_1}^{\mu_1}
x_{i_2}^{\mu_2}\cdots x_{i_{\ell(\mu)}}^{\mu_{\ell(\mu)}}$, each
monomial corresponds uniquely to the labelling of ${\cal I}$ where the
$j$th interval is labelled by $i_j$, and so occurs exactly once in
$\sum x^{{\cal I} }$ as required.
\end{proof}

The \emph{augmented monomial symmetric function} $\tilde{m}_\mu$ is
defined by 
  $$ \tilde{m}_\mu = m_1(\mu)!m_2(\mu)!\cdot m_\mu, $$
where $m_\mu$ denotes a monomial symmetric function. 

\begin{lemma}\label{l:monomials}
$c_\mu \tilde{m}_\mu  = (-1)^{z(\lambda)} \sum
(-1)^{c( {\cal I})}  x^{
  {\cal I}}$, where the sum is over all labelled interval sets of
$\lambda$ of type $\mu$ with no label repeated.  
\end{lemma}

Note that we have already demonstrated this result in Example
\ref{e:powersums}.  Indeed for that interval set when the labels were
all different we saw 
that  $\sum  x^{  {\cal I}} = 2m_{633}$.

\begin{proof}
Fix a specific interval set ${\cal I}$ of type $\mu$.  We need to show
$\tilde{m}_\mu = \sum x^{ {\cal I}}$ where the sum is over all
labellings with no label repeated.  As before we can order the
intervals and say that the $i$th interval is labelled by $\alpha_i$.
Note that if $\mu_j = \mu_{j+1}$, the two labellings $(\alpha_1,
\alpha_2, \ldots, \alpha_j, \alpha_{j+1}, \ldots)$ and $(\alpha_1,
\alpha_2, \ldots, \alpha_{j+1}, \alpha_{j}, \ldots)$ both produce the
same term $x^{ {\cal I}} = x_{\alpha_1}^{\mu_1}
x_{\alpha_2}^{\mu_2}\cdots$.  So we have $\sum_{(\alpha_1, \alpha_2,
\ldots)} x^{ {\cal I}} =\sum_{(\beta_1, \beta_2, \ldots)} m_1(\mu)!
m_2(\mu)! \cdots x^{ {\cal I}} $ where we impose the condition that if
$\mu_j = \mu_{j+1}$, then $\beta_j < \beta_{j+1}$.  Recall that
$m_\mu= \sum_{(\beta_1, \beta_2, \ldots)} x^{ {\cal I}}$ by
definition.  So we have $\tilde{m}_\mu = \sum x^{ {\cal
I}}_{(\alpha_1, \alpha_2, \ldots)}$ as required.
\end{proof}

\begin{theorem}\label{th:identity}
For each shape $\lambda$, write the bottom Schur function
$\hat{s}_\lambda = \sum_\mu c_\mu \tilde{p}_\mu$.  Then
$\sum_\mu c_\mu p_\mu = \sum_\mu c_\mu \tilde{m}_\mu$.  
\end{theorem}

\begin{example}
For $\lambda = (4,4,4)$ we have 
$$\hat{s}_\lambda = 
-\tilde{p}_{642}
+\tilde{p}_{633}
+\tilde{p}_{552}
-2 \tilde{p}_{543}
+\tilde{p}_{444}.$$
  So our result states that
$$
-p_{642}+p_{633}+p_{552}-2p_{543}+p_{444} = 
-m_{642}+2m_{633}+2m_{552}-2m_{543}+6m_{444}.
$$

\end{example}

\begin{proof}
From Lemma \ref{l:powersums} we have $c_\mu p_\mu  =
(-1)^{z(\lambda)} \sum (-1)^{c( {\cal I})} x^{ 
  {\cal I}}$, where the sum is over all labelled interval sets of type
$\mu$.  But by Lemma \ref{l:indistinct} 
  $$ (-1)^{z(\lambda)} \sum (-1)^{c( {\cal I})} x^{ 
  {\cal I}} = 0 $$ 
if we sum over all labelled interval sets with a repeated label, while 
by Lemma \ref{l:monomials} we have
 $$ (-1)^{z(\lambda)}  \sum (-1)^{c( {\cal I})} x^{
  {\cal I}} =  c_\mu m_1(\mu)! m_2(\mu)! \cdots m_\mu $$ 
if we sum over all labelled interval sets with no label repeated.  So
$\sum_\mu c_\mu p_\mu = \sum_\mu c_\mu \tilde{m}_\mu$.  
\end{proof}

Let $\Gamma$ denote the space of all symmetric functions $f$ with
rational coefficients such that if $f=\sum_\mu c_\mu \tilde{p}_\mu$,
then 
  $$ \sum_\mu c_\mu p_\mu = \sum_\mu c_\mu\tilde{m}_\mu. $$
Theorem~\ref{th:identity} shows that
$\hat{s}_\lambda\in\Gamma$. However, there are elements of $\Gamma$
that are not linear combinations of $\hat{s}_\lambda$'s, such as $f =
\tilde{p}_{511}
-3\tilde{p}_{421}+\tilde{p}_{331}+\tilde{p}_{322}$. Let $R$ denote the
transition matrix from monomial symmetric functions to the power sums,
i.e., 
   $$ p_\lambda = \sum_\mu R_{\lambda\mu}m_\mu. $$
Let $D$ denote the diagonal matrix whose diagonal coincides with that
of $R$, so
  $$ D_{\lambda\lambda} = R_{\lambda\lambda} = \prod_{i\geq 1}
     m_i(\lambda)!. $$
It is easy to see that $\Gamma=\ker(R-D)$, the kernel (or null space) of
$R-D$. Let $\Gamma_n$ denote the elements of $\Gamma$ that are
homogeneous of degree $n$, and let $\gamma_n=\dim \Gamma_n$. We have
computed that 
  $$ (\gamma_1,\gamma_2, \dots) = (1,1,1,2,2,3,4,5,7,9,11,15,19,24,
        \dots). $$
Compare this with the dimension $\beta_n$ of the space spanned by the
bottom Schur functions of degree $n$, given in general by
Corollary~\ref{cor:botdim} and for $n\leq 27$ just below this
corollary. In particular, the least $n$ for which $\beta_n<\gamma_n$
is $n=7$. We don't have a conjecture for the value of $\gamma_n$.

\appendix



\begin{thebibliography}{99}

\bibitem{Na98} M. Nazarov and V. Tarasov, On irreducibility of
  tensor products of Yangian modules, \emph{Internat. Math. Research
    Notices} (1998), 125-150.


\bibitem{St99} R. P. Stanley, ``Enumerative Combinatorics,'' Vol. 2,
  Cambridge Univ. Press, New York/Cambridge, UK, 1999.

\bibitem{St02} R. P. Stanley, The Rank and Minimal Border Strip
Decompositions of a Skew Partition, \emph{J. Combin. Theory (A)} {\bf
100} (2002), 349-375.

\bibitem{St01} J. R. Stembridge, The SF Package for Maple, Version
  2.3, 22 July 2001,
  http://www.math.lsa.umich.edu/$\sim$jrs/maple.html\#SF

\end{thebibliography}
\end{document}